\documentclass[12pt]{amsart}
\usepackage{graphicx}
\usepackage{amssymb}
\usepackage{amsfonts}
\usepackage{amsmath}
\usepackage{array}
\usepackage{rotating}
\usepackage{bm}
\usepackage[matrix,frame]{xypic}
\usepackage{MnSymbol}
\newtheorem{example}{Example}[section]

\newtheorem{theorem}[example]{Theorem}

\newtheorem{definition}[example]{Definition}
\newtheorem{proposition}[example]{Proposition}

\def\Proof{\noindent \it Proof -- \rm}
\def\qed{\hspace{3.5mm} \hfill \vbox{\hrule height 3pt depth 2 pt width 2mm}
\bigskip}

\def\WQSym{{\bf WQSym}}

\def\<{\langle}
\def\>{\rangle}

\def\SG{{\mathfrak S}}

\def\sinv{{\rm sinv}}

\def\std{{\rm std}}

\def\LL{{\mathfrak P}}

\title[]%
{On composition polynomials}
\author[J.-C.~Novelli and J.-Y.~Thibon]%
{Jean-Christophe Novelli and Jean-Yves Thibon}

\address[]{[Novelli, Thibon] Laboratoire d'informatique Gaspard-Monge\\
Universit\'e Paris-Est Marne-la-Vall\'ee \\
5, Boulevard Descartes \\ Champs-sur-Marne \\
77454 Marne-la-Vall\'ee cedex 2 \\
France}
\email[Jean-Christophe Novelli]{novelli@univ-mlv.fr}
\email[Jean-Yves Thibon]{jyt@univ-mlv.fr} 
\date{\today}


\date{\today}

\begin{document}

\begin{abstract}
We provide a combinatorial interpretation of the reduced composition polynomials
of Ardila and Doker [Adv. Appl. Math. {\bf 50} (2013), 607], and relate them
to the $(1-q)$-transform of noncommutative symmetric functions. 
\end{abstract}

\maketitle

\section{Composition polynomials and noncommutative symmetric functions}

Let $I=(i_1,\ldots,i_r)$ be a composition of $n$.
The {\it composition polynomial} $g_I(q)$ is defined in \cite{AD} as
\begin{equation}
g_I(q)=\int_q^{1}\int_q^{t_r}\cdots\int_q^{t_2}t_1^{i_1-1}t_2^{i_2-1}\cdots t_r^{i_r-1}dt_1dt_2\cdots dt_r.
\end{equation}
This polynomial is divisible by $(1-q)^r$, and 
$f_I(q)=g_I(q)/(1-q)^r$ is called the {\it reduced composition polynomial}.
It is proved in \cite{AD} that its coefficients are nonnegative. Our aim
is to give a combinatorial interpretation of the coefficients of $n!f_I(q)$,
which turn out to be integers.
Actually, these polynomials have a natural interpretation in
the context of noncommutative symmetric functions \cite{NCSF1,NCSF2}.

Let $S_n(A)$ denote the noncommutative complete symmetric functions of an
alphabet $A$, and $\Psi_n(A)$ the noncommutative power-sums of the first
kind.

Their generating series
\begin{equation}
\sigma_t(A) := \sum_{n\ge 0}S_n(A)t^n\quad\text{and}\quad \psi_t(A) :=\sum_{n\ge 1}\Psi_nt^{n-1}
\end{equation}
are related by the differential equation
\begin{equation}
\frac{d}{dt}\sigma_t(A)= \sigma_t(A)\psi_t(A)
\end{equation}
with the initial condition $\sigma_0(A)=1$. This can be recast as an integral equation
\begin{equation}
\sigma_t(A)= 1+\int_0^t\sigma_u(A)\psi_u(A)du
\end{equation}
which may be solved by iterated integrals\footnote{The reader unfamiliar with noncommutative symmetric function
can assume that $\sigma_t$ is a generic noncommutative power series and take this as a definition of $\psi_t$.}.

In \cite{NCSF2}, an analogue of the classical $(1-q)$-transform of ordinary
symmetric functions (sending the power sums $p_n$ to $(1-q^n)p_n$) is defined:
\begin{equation}
\sigma_t((1-q)A) := \lambda_{-qt}(A)\sigma_t(A)\quad (\text{where}\ \lambda_t(A)=\sigma_t(A)^{-1}).
\end{equation}

\begin{proposition}
The coefficient of
$\Psi^I$ in  $S_n((1-q)A)$ is $g_I(q)$.
\end{proposition}

\Proof The series $X(t)= \sigma_1((t-q)A):=\lambda_{-q}(A)\sigma_t(A)$ 
satisfies the differential equation
\begin{equation}
\frac{dX}{dt} = \lambda_{-q}(A)\frac{d}{dt}\sigma_t(A)=\lambda_{-q}(A)\sigma_t(A)\psi_t(A)
=X(t)\psi_t(A)
\end{equation}
with the initial condition
\begin{equation}
X(q)=1,
\end{equation}
so that
\begin{equation}
X(t)=1+\int_q^t X(t_1)\psi_{t_1}dt_1
\end{equation}
whose solution is
\begin{equation}
X(t) = 1 + \sum_{r\ge 1}
\int_q^t\int_q^{t_1}\cdots\int_q^{t_{r-1}}\psi_{t_r}\cdots\psi_{t_2}\psi_{t_1}dt_r\cdots dt_2dt_1.
\end{equation}
\qed

\begin{example}\label{exS4}{\rm For $n=4$,
\begin{equation}\label{eq:S4}
\begin{split}
4!S_4((1-q)A) = & 
\ 6\, \left( 1-q \right) \left( q+1 \right)  \left( {q}^{2}+1 \right) { \Psi}^{{4}}
+2\, \left( 3\,{q}^{2}+2\,q+1 \right)  \left( 1-q \right) ^{2}{ \Psi}^{{31}}\\
& +3\, \left( 1-q \right) ^{2} \left( q+1 \right) ^{2}{ \Psi}^{{22}}
+ \left( 3\,q+1 \right)  \left( 1-q \right) ^{3}{ \Psi}^{{211}}\\
& +2\, \left( {q}^{2}+2\,q+3 \right)  \left( 1-q \right) ^{2}{ \Psi}^{{13}}
+2\, \left( q+1 \right)  \left( 1-q \right) ^{3}{ \Psi}^{{121}}\\
& + \left( q+3 \right)  \left( 1-q \right) ^{3}{ \Psi}^{{112}}
+ \left( 1-q \right) ^{4} { \Psi}^{{1111}}.
\end{split}
\end{equation}
}
\end{example}

\begin{definition}
For a composition $I=(i_1,\ldots,i_r)$ of $n$, we set
\begin{equation}
 P_I(q)=n!f_I(q)=n!\frac{g_I(q)}{(1-q)^r}.
\end{equation}
\end{definition}

\begin{example}{\rm
From Example~\ref{exS4}, we have
\begin{equation}
\label{ex-pp}
\begin{array}{lll}
P_{4}   &=& 6q^3+6q^2+6q+6,  \\
P_{31}  &=& 6q^2+4q+2     ,  \\
P_{22}  &=& 3q^2+6q+3     ,  \\
P_{211} &=& 3q+1          ,  \\
P_{13}  &=& 2q^2+4q+6     ,  \\
P_{121} &=& 2q+2          ,  \\
P_{112} &=&  q+3          ,  \\
P_{1111}&=& 1             .  \\
\end{array}
\end{equation}
}
\end{example}

Note that the sequence of the sums of all $P_I$ of a given weight at $q=1$
begins with
\begin{equation}
1,1,3,13,73
\end{equation}
which suggests a combinatorial interpretation using sets of lists.

\section{Combinatorial interpretation of the reduced polynomials}

A \emph{set of lists} of size $n$ is a set partition $\pi$
of $[n]$ endowed with a total order inside each block \cite[A000262]{Sloane}.
Since this object is intermediate between a segmented permutation
and a set partition, 
we propose to call it a {\it permutation partition}, or for short, a {\it permutition}.

Let $\LL_n$ be the set of permutitions of size $n$.
We define a statistic $\sinv$ (special inversions) on $\LL_n$ as follows:
$\sinv(\pi)$ is the number of letters in $\pi$ greater than the last letter of
their block.

We also introduce an order among the blocks, according to the value of the
last letter of each block, and define the composition $c(\pi)$ as the sequence
of lengths of the blocks ordered in this way.

For example, with $\pi=\{53, 4612,978\}$, $\sinv(\pi)=1+2+1=4$,
and the canonical ordering of $\pi$ is the segmented permutation
$\sigma=(4612|53|978)$, so that $c(\pi)=(4,2,3)$. 

\begin{theorem}
The reduced composition polynomial $P_I$ is the generating
function of $\sinv$ on permutitions of composition $I$:
\begin{equation}\label{stat}
P_I(q)=\sum_{c(\pi)=I}q^{\sinv(\pi)}.
\end{equation}
\end{theorem}

\begin{example}{\rm
With $I=(1,3)$, one gets 12 permutitions of shape $I$:
$$(1|234), (1|324), (1|243), (1|423), (1|342), (1|432), $$
$$ (2|134), (2|314), (2|143), (2|413),$$
$$ (3|124), (3|214)$$
whose numbers of special inversions are respectively $0,0,1,1,2,2,0,0,1,1,0,0$,
yielding the polynomial $6+4q+2q^2$ which is indeed $P_{13}$, see
Example~\ref{ex-pp}.
}
\end{example}

\bigskip
\Proof 
Ardila and Doker \cite{AD} give the induction formula
\begin{equation}
\frac1{i_1}g_{(i_1+i_2,i_3,\ldots,i_r)}(q)
= g_I(q)-\frac{q^{i_1}}{i_1}g_{(i_2,\ldots,i_r)}(q)
\end{equation}
which translates into
\begin{equation}
\label{AD-rec}
qi_1P_I(q) +
P_{(i_1+i_2,i_3,\ldots,i_r)}(q) 
=
i_1P_I(q) + \frac{n!}{(n-i_1)!}q^{i_1}P_{(i_2,\ldots,i_r)}(q),
\end{equation}
which, together with the initial conditions $P_{1^n}(q)=1$, completely
determine all the $P_I$.
To prove that the polynomials defined by \eqref{stat} satisfy this induction,
we interpret both sides of \eqref{AD-rec}
as $q$-enumerations of disjoint unions of sets.

\medskip
Let us first describe the combinatorial interpretation of each of the four terms.

First, both terms $i_1 P_I(q)$ are seen as the set $E_1$ of pairs $(\tau,x)$
where $\tau$ satisfies $c(\tau)=I$ and $x$ is an element of 
$\{\tau_1,\dots,\tau_{i_1}\}$.

Following Ardila-Doker,  denote by $I^1$ the composition
$(i_1+i_2,i_3,\ldots,i_r)$.
The term $P_{I^1}$ is interpreted as the set $E_2$ of permutitions $\tau$
satisfying $c(\tau)=I^1$.

Finally, the term $\frac{n!}{(n-i_1)!}P_{(i_2,\ldots,i_r)}$ is translated as
the set $E_3$ of pairs $(\tau,w)$ where $\tau$ satisfies
$c(\tau)=(i_2,\ldots,i_r)$ and where $w=w_1\dots w_{i_1}$ is a word of size
$i_1$ with distinct values between $1$ and $n$.

\medskip
We shall now proceed as follows: first define a map from $E_1$ to  $E_1\cup E_3$, 
characterize its image set, show that it is
bijective and show that it behaves as desired with respect to the $q$-enumeration
constraint. Then, define a map from $E_2$ to  $E_1\cup E_3$, show the same
properties as before and also show that the images of both partition
the set $E_1\cup E_3$.

\medskip
Let us first define our map $\phi_1$ from $E_1$.
Consider an element $(\tau,x)$ of $E_1$
and  set $S=\{\tau_1,\dots,\tau_{i_1}\}$.
\begin{itemize}
\item If $\tau_{i_1}$ is the smallest value of $S$,
define $w$ as the word obtained from $\tau_1\dots\tau_{i_1}$ by exchanging
$x$ with $\tau_{i_1}$. Then
\begin{equation}
\phi_1((\tau,x))
= (\std(\tau_{i_1+1}\dots \tau_{n}), w),
\end{equation}
where $\std$ denotes the usual standardization of words over an ordered alphabet.
\item Otherwise, let $m$ be the maximal value of $S$ strictly smaller than
$\tau_{i_1}$. Then let $\tau'$ be obtained from $\tau$ by exchanging
$\tau_{i_1}$ with $m$ and
define
\begin{equation}
\phi_1((\tau,x)) = (\tau',x).
\end{equation}
\end{itemize}
By construction, 
in both cases, each element in the image of
$\phi_1$ has a unique preimage,
 so that
$\phi_1$ is a bijection from $E_1$ to its image set. Indeed,
one can retrieve $\tau$ in the first case, since the values
$\tau_{i_1+1}\dots \tau_{n}$ are obtained from their standardized word
by the unique increasing morphism from $[n-i_1]$ to the set of values 
in $[n]$
not occuring in $w$. 
Then, $\tau_{1}\dots\tau_{i_1}$ is obtained from
$w$ by exchanging the smallest value of $w$ with its last one, and  $x$ is
$w_{i_1}$. The second case is straightforward, as $x$ plays no r\^ole in
the definition of $\tau'$, from which $\tau$ can be easily recovered.

Let us now describe the image set $\phi_1(E_1)$.
The intersection of $E_1$ and $\phi_1(E_1)$ is the set of
elements where $x$ can take any value and such that $\tau_{i_1}$ is not
the greatest value of  $\{\tau_1,\dots,\tau_{i_1}\}$ smaller than
$\tau_{i_1+i_2}$. 
Indeed, such a pair clearly has a preimage as described above and the other
elements of $E$ do not. The intersection of $E_3$ after applying the
increasing morphism to the standard word and $\phi_1(E_1)$ is the set of
elements where the smallest value of $w$ is smaller than $\tau_{i_1+i_2}$.
Indeed, the condition on the smallest value of $w$ is necessary
and sufficient to have a permutition as preimage.

Finally, with respect to the $q$-enumeration according to the number of
special inversions of the permutitions on both sides, $\phi_1$ behaves as
expected.
Any element $(\tau,x)\in E_1$ sent to $(\pi,x)\in E_1$ satisfies 
$\sinv(\tau)=\sinv(\pi)-1$, since there will be one more inversion in $\pi$
than in $\tau$ among their first $i_1$ values, that is, before the first bar.
Now, for any $(\tau,x)\in E_1$ sent to $(\pi,w)\in E_3$, 
$\sinv(\tau)=\sinv(\pi)+i_1-1$ since there will be $i_1-1$ less inversions  in
$\pi$ than in $\tau$: indeed, $\tau_{i_1}$ is the smallest value before its
first bar, so that there are $i_1-1$ inversions in $\tau$ before its first bar.

\bigskip
%
Let us now define a second map $\phi_2$ from  $E_2=\{\tau| c(\tau)=I^1\}$
to $E_1\cup E_3$.
\begin{itemize}
\item If $\tau_\alpha>\tau_{i_1+i_2}$ for all $\alpha\leq i_1$,
then 
\begin{equation}
\phi_2(\tau)
=(\std(\tau_{i_1+1}\dots\tau_{i_1+i_2}|\dots),
  \tau_{1}\dots\tau_{i_1}).
\end{equation}
\item Otherwise, let $k$ be the greatest element of
$\{\tau_1,\dots,\tau_{i_1}\}$ smaller than $\tau_{i_1+i_2}$
and let $\tau'$ be obtained from ${\tau_{1}\dots\tau_{i_1}}$
by exchanging $k$ with $\tau_{i_1}$, and define
\begin{equation}
\phi_2(\tau)
=((\tau' | \tau_{i_1+1}\dots\tau_{i_1+i_2}|\dots),
  \tau_{i_1}).
\end{equation}
\end{itemize}

As with $\phi_1$, $\phi_2$ is a bijection from $E_2$ to its image set.

Let us now describe the image set  $\phi_2(E_2)$.
The intersection of $E_1$ and  $\phi_2(E_2)$ is the set of elements such that
$x$ can take any value and such that $\tau_{i_1}$ is the greatest element of
$\{\tau_1,\dots,\tau_{i_1}\}$ smaller than $\tau_{i_1+i_2}$.
The intersection of $E_3$ and $\phi_2(E_2)$ is the set of
elements such that the smallest value of $w$ is greater than $\tau_{i_1+i_2}$.
Hence,  both sets constitute the  complement of the image of $\phi_1$, so that
the map $\phi$ defined by $\phi|_{E_1}=\phi_1$ and $\phi|_{E_2}=\phi_2$ 
is a bijection between $E_1\cup E_2$
and $E_1\cup E_3$.

Again, $\phi_2$ behaves as expected with respect to the statistic $\sinv$.
Any element $\tau\in E_2$ sent to $(\pi,x)\in E_1$ satisfies 
$\sinv(\tau)=\sinv(\pi)$, since the inversions between $\tau_k$ with
$k\leq i_1$ and $\tau_{i_1+i_2}$ become all inversions between $\pi_i$ with
$i\leq i_1$ and $\pi_{i_1}$.
Moreover, for any $\tau\in E_2$ sent to $(\pi,w)\in E_3$, 
$\sinv(\tau)=\sinv(\pi)+i_1$, since there are inversions between all
$\tau_k$ with $k\leq i_1$ and $\tau_{i_1+i_2}$ that disappear in $\pi$.
Therefore, $\phi$ is a bijection between $E_1\sqcup E_2$ and
$E_1\sqcup E_3$, behaving as claimed with respect to $\sinv$.
\qed

\subsection{Examples}

\subsubsection{}
Let us first give four examples 
of $\phi_1$ and $\phi_2$ and their inverse maps, covering all possible cases.

We have
\begin{equation}
\begin{split}
\phi_1((361|74|258),3) &= ((42|135),163) \\
\phi_1((163|74|258),1) &= ((361|74|258),1) \\
\phi_2((26371|458))    &= ((41|235),263) \\
\phi_2((36174|258))    &= ((163|74|258),1)
\end{split}
\end{equation}
Now, from $((42|135),163)$, one easily recovers the suffix of the permutition $\tau$ of its preimage.
It is $(74|258)$, from which one can see that $163$ are not all
greater than $4$. Hence, it belongs to the image of $\phi_1$, and then easily
gives back the first permutition of the example. 

In the case of
$((41|235),263)$, one recovers as suffix the expression $(71|458)$, and since
all values of $263$ are greater than $1$, it belongs to the image of $\phi_2$.
Finally, on the other two examples $((361|74|258),1)$ and $((163|74|258),1)$,
since $3$ is the greatest possible value of the prefix considering that the
second expression between bars ends with a $4$, these elements respectively
belong to the image of $\phi_1$ and $\phi_2$.

\subsubsection{}
Let us now illustrate the complete bijection on the example of the composition $I=(2,1,1)$.
The following table shows the images of $E_1=\{(\tau,x)\}$ satisfying
$c(\tau)=(2,1,1)$ and $x\in\{1,2\}$ with the corresponding $q$-statistics on the
right:
 
\begin{equation}
\label{ex-phi1}
\begin{split}
(12|3|4, 1) &\mapsto (21|3|4,1) \qquad\qquad (1,q) \\
(12|3|4, 2) &\mapsto (21|3|4,2) \qquad\qquad (1,q) \\
(21|3|4, 1) &\mapsto (1|2,21) \qquad\qquad (q,1) \\
(21|3|4, 2) &\mapsto (1|2,12) \qquad\qquad (q,1) \\
(31|3|4, 1) &\mapsto (1|2,31) \qquad\qquad (q,1) \\
(31|3|4, 3) &\mapsto (1|2,13) \qquad\qquad (q,1) \\
(41|3|4, 1) &\mapsto (1|2,41) \qquad\qquad (q,1) \\
(41|3|4, 4) &\mapsto (1|2,14) \qquad\qquad (q,1) \\
\end{split}
\end{equation}

The next table shows the images of $E_2=\{\tau\}$ satisfying $c(\tau)=(3,1)$
with the corresponding $q$-statistics on the right:
\begin{equation}
\label{ex-phi2}
\begin{split}
(123|4) &\mapsto (12|3|4,2) \qquad\qquad (1,1) \\
(213|4) &\mapsto (12|3|4,1) \qquad\qquad (1,1) \\
(132|4) &\mapsto (31|2|4,3) \qquad\qquad (q,q) \\
(312|4) &\mapsto (31|2|4,1) \qquad\qquad (q,q) \\
(142|3) &\mapsto (41|2|3,4) \qquad\qquad (q,q) \\
(412|3) &\mapsto (41|2|3,1) \qquad\qquad (q,q) \\
(231|4) &\mapsto (1|2, 23) \qquad\qquad (q^2,1) \\
(321|4) &\mapsto (1|2, 32) \qquad\qquad (q^2,1) \\
(241|3) &\mapsto (1|2, 24) \qquad\qquad (q^2,1) \\
(421|3) &\mapsto (1|2, 42) \qquad\qquad (q^2,1) \\
(341|2) &\mapsto (1|2, 34) \qquad\qquad (q^2,1) \\
(431|2) &\mapsto (1|2, 43) \qquad\qquad (q^2,1) \\
\end{split}
\end{equation}

Now, summing up the $q$-statistics on both side of~\eqref{ex-phi1}
and~\eqref{ex-phi2}, one gets that its left-hand side is
\begin{equation}
q(6q+2) + (6q^2+4q+2) = 12q^2+6q+2,
\end{equation}
and its right-hand side is
\begin{equation}
6q+2 + q^2\times 12,
\end{equation}
so that they coincide.

\section{Miscellaneous comments}

\subsection{}Noncommutative symmetric functions of degree $n$ can be
intepreted as elements of the descent algebra $\Sigma_n$ of $\SG_n$.
In this context, $\Psi_n$ is the Dynkin symmetrizer (iterated
bracketing): as a linear combination of permutations, 
\begin{equation} 
\Psi_n = [\cdots [[1,2],3], \cdots,n]
\end{equation}
and $S_n((1-q)A$ is the iterated $q$-bracketing
\begin{equation} 
S_n((1-q)A)= (1-q) [\cdots [[1,2]_q,3]_q, \cdots,n]_q.
\end{equation}
By \cite[Lemma 5.11]{NCSF2}, writing
\begin{equation}
S_n((1-q)A) = S_n((1-q)A)*S_n(A)
\end{equation}
(internal product) and inserting the expansion
\begin{equation}
S_n(A) = \sum_{r=1}^n\sum_{I\vDash n\atop \ell(I)=r}\frac{\Psi^I}{i_1(i_1+i_2)\cdots(i_1+\cdots i_r)}
\end{equation}
we obtain
\begin{equation}
S_n((1-q)A) = \sum_{r=1}^n\sum_{I\vDash n\atop \ell(I)=r}(1-q^{i_1})
\frac{[\cdots [[\Psi_{i_1},\Psi_{i_2}]_{q^{i_2}},\Psi_{i_3}]_{q^{i_3}},\cdots,\Psi_{i_r}]_{q^{i_r}}} 
{i_1(i_1+i_2)\cdots(i_1+\cdots i_r)}
\end{equation}

\subsection{}
There exists a combinatorial Hopf algebra based on permutitions \cite{BCLM}.
It plays with respect to the Hopf algebra of set partitions (symmetric
functions in noncommuting variables) a r\^ole symmetrical to that of $\WQSym$
(quasi-symmetric functions in noncommuting variables).

The number of permutitions of length $n$ is equal to the number of stalactic
classes of parking functions of the same length, \cite{HNT}, on which a
combinatorial Hopf algebra structure can also be defined. 
There is also a bijection between permutitions and non-crossing set
compositions \cite{Cal}. The bijection goes through Dyck paths with labelled
peaks, which are easily identified with the canonical representatives of
stalactic classes of parking functions, obtained by permuting among themselves
the blocks of identical letters in a nondecreasing parking function. 
Composing these bijections, on can obtain a bijection between permutitions and
stalactic classes of parking functions.



\footnotesize
\begin{thebibliography}{aa}
\bibitem{AD}{\sc F. Ardila} and {\sc J. Doker}, {\it
Lifted generalized permutahedra and composition polynomials},
Adv. Appl. Math. {\bf 50} (2013), 607--633.
%
\bibitem{BCLM}{\sc J.-P. Bultel, A. Chouria, J.-G. Luque}, and {\sc O. Mallet},
{\it Word symmetric functions and the Redfield-P\'olya theorem},
25th International Conference on Formal Power Series and Algebraic
Combinatorics (FPSAC 2013), 563--574,
Discrete Math. Theor. Comput. Sci. Proc., AS, Assoc. Discrete Math. Theor.
Comput. Sci., Nancy, 2013.
%

\bibitem{Cal}{\sc D. Callan}, {\it
Sets, lists and noncrossing partitions}, J. Integer Seq. {\bf 11} (2008), no. 1, Article 08.1.3, 7 pp.
%
\bibitem{HNT}{\sc F. Hivert, J.-C. Novelli}, and {\sc J.-Y. Thibon},
{\it Commutative combinatorial Hopf algebras},
J. Algebraic Combin. {\bf 28} (2008), no. 1, 65--95.
%

\bibitem{NCSF1} {\sc I.~M. Gelfand, D. Krob, A. Lascoux, B. Leclerc,
V.~S. Retakh}, and {\sc J.-Y. Thibon}.
{\it Noncommutative symmetric functions},
Adv. Math. {\bf 112}, 1995, 218--348.
%
%
\bibitem{NCSF2} {\sc Krob D., Leclerc B.}, and {\sc Thibon J.-Y.},
{\it Noncommutative symmetric functions II: Trans\-for\-ma\-tions of
alphabets}, Intern. J. Alg. Comput.,
{\bf 7}, (2), (1997), 181--264.
%
%
\bibitem{Sloane}{\it The On-Line Encyclopedia of Integer Sequences}, 
published electronically at http://oeis.org, 2010. 
%

\end{thebibliography}
\end{document}